\documentclass[11pt]{article}
\usepackage{amsmath,amsfonts,latexsym,amssymb,amscd, graphicx,pictexwd}
\renewcommand{\P}{{\mathbb P}}
\newcommand{\E}{{\mathbb E}}

\newcommand{\R}{{\mathbb R}}

\newcommand{\ZZ}{{\mathbf Z}}

\newcommand{\cc}{{\mathbf c}}

\newcommand{\cA}{{\mathcal A}}

\newcommand{\cB}{{\mathcal B}}

\newcommand{\Ai}{{\rm Ai}}
\newcommand{\eps}{{\varepsilon }}

\newcommand{\la}{{\lambda}}

\newtheorem{thm}{Theorem}

\newtheorem{lem}{Lemma}

\newtheorem{rem}{Remark}
\numberwithin{equation}{section}

\begin{document}
\title{On the local fluctuations of last-passage percolation models}
\author{Eric Cator and Leandro P. R. Pimentel}
\date{\today}
\maketitle

\begin{abstract}  
Using the fact that the Airy process describes the limiting fluctuations of the Hammersley last-passage percolation model, we prove that it behaves locally like a Brownian motion. Our method is quite straightforward, and it is based on a certain monotonicity and good control over the equilibrium measures of the Hammersley model (local comparison).
\end{abstract}

\section{Introduction and results}

In recent years there has been a lot of research on the Airy process and related processes such as the Airy sheet \cite{CQ}. Most papers use analytic methods and exact formulas given by Fredholm determinants to prove properties of these processes, but some papers use the fact that these processes are limiting processes of last-passage percolation models or random polymer models, and they use properties of these well studied models to prove the corresponding property of the limiting process. A nice example of these different approaches can be found in two recent papers, one by H\"agg \cite{H} and one by Corwin and Hammond \cite{CH}. H\"agg proves in his paper that the Airy process behaves locally like a Brownian motion, at least in terms of convergence of finite dimensional distributions. He uses the Fredholm determinant description of the Airy process to obtain his result. Corwin and Hammond on the other hand, use the fact that the Airy line process, of which the top line corresponds to the Airy process, can be seen as a limit of Brownian bridges, conditioned to be non-intersecting. They show that a particular resampling procedure, which they call the Brownian Gibbs property, holds for the system of Brownian bridges and also in the limit for the Airy line process. As a consequence, it follows that the local Brownian behavior of the Airy process holds in a stronger functional limit sense.  Our paper will prove the same theorem, also using the fact that the Airy process is a limiting process, but in a much more direct way: we will consider the Hammersley last-passage percolation model, and show that we can control local fluctuations of this model by precisely chosen equilibrium versions of this model, which are simply Poisson processes. Then we show that in the limit this control suffices to prove the local Brownian motion behavior of the Airy process, as well as tightness of the models approaching the Airy process. We also extend the control of local fluctuations of the Hammersley process to scales smaller than the typical cube-root scale.

Our method is quite straightforward, yet rather powerful, mainly because we have a certain monotonicity and good control over the equilibrium measures. In fact, we think that we can extend our result to the more illustrious Airy sheet, the two dimensional version of the Airy process. We address the reader to \cite{CQ}, for a description of this process in terms of the renormalization fixed point of the KPZ universality class. However, here we run in to much more technical problems, and this will still require a lot more work, beyond the scope of this paper.

We will continue the introduction by developing notation, introducing all relevant processes and stating the three main theorems. In Section 2 we introduce the local comparison technique and in each of the following three sections one theorem is proved.

\subsection{The Hammersley Last Passage Percolation model}

The Hammersley last-passage percolation model \cite{AD} is constructed from a two-dimensional homogeneous Poisson point process of intensity $1$. Denote $[x]_t:=(x,t)\in\R^2$ and call a sequence $[x_1]_{t_1},[x_2]_{t_2},\dots,[x_k]_{t_k}$ of planar points increasing if $x_j<x_{j+1}$ and $t_j<t_{j+1}$ for all $j=1,\dots,k-1$. The last-passage time $L([x]_s,[y]_t)$ between  $[x]_s<[y]_t$ is the maximal number of Poisson points among all increasing sequences of Poisson points lying in the rectangle $(x, y ]\times(s, t]$. Denote $L[x]_t:=L([0]_0,[x]_t)$ and define $\cA_n$ by 
$$u\in\R\,\mapsto\,\cA_n(u):=\frac{L[n+2un^{2/3}]_n-(2n+2un^{2/3})+u^2n^{1/3}}{n^{1/3}}\,.$$
Pr\"ahofer and Spohn \cite{PS} proved that 
\begin{equation}\label{eq:Airy}
\lim_{n\to\infty}\cA_n(\cdot)\stackrel{dist.}{=}\cA(\cdot)\,,
\end{equation}
in the sense of finite dimensional distributions, where $\cA\equiv(\cA(u))_{u\in\R}$ is the so-called Airy process. This process is a one-dimensional stationary process with continuous paths and  finite dimensional distributions given by a Fredholm determinant \cite{J}:
$$\P\left(\cA(u_1)\leq \xi_1,\dots,\cA(u_m)\leq \xi_m\right):=\det\left(I-f^{1/2}Af^{1/2}\right)_{L^2\left(\{u_1,\dots,u_m\}\times\R\right)}\,.$$
The function $A$ denotes the extended Airy kernel, which is defined as
$$A_{s,t}(x,y):=\left\{\begin{array}{ll}\int_0^\infty e^{-z(s-t)}\Ai(x+z)\Ai(y+z)dz\,,& \mbox{ if } s\geq t\,,\\&\\
-\int_{-\infty}^0 e^{-z(t-s)}\Ai(x+z)\Ai(y+z)dz\,,& \mbox{ if } s< t\,,\end{array}\right.$$
where $\Ai$ is the Airy function, and for $\xi_1,\dots,\xi_m\in\R$ and $u_1<\dots<u_m$ in $\R$, 
\begin{eqnarray*}
f\,:\,\{u_1,\dots,u_m\}\times\R&\to&\R\\
(u_i,x)&\mapsto&1_{(\xi_i,\infty)}(x)\,.
\end{eqnarray*}

The main contribution of this paper is the development of a local comparison technique to study the local fluctuations of last-passage times and its scaling limit. The ideas parallel the work of Cator and Groeneboom \cite{CG1,CG2}, where they studied local convergence to equilibrium and the cube-root asymptotic behavior of $L$. This technique consists of bounding from below and from above the local differences of $L$ by the local differences of the equilibrium regime (Lemma \ref{lem:LocalComparison}), with suitable  parameters that will allow us to handle the local fluctuations in the right scale (Lemma \ref{lem:ExitControl}). For the Hammersley model the equilibrium regime is given by a Poisson process. We have strong indications that the technique can be applied to a broad class of models, as soon as one has Gaussian fluctuations for the equilibrium regime. Although this is a very natural assumption, one can only check that for a few models. As a first application, we will prove tightness of $\cA_n$.   

\begin{thm}\label{thm:Tight}
The collection $\{\cA_n\}$ is tight in the space of cadlag functions on $[a,b]$. Furthermore, any weak limit of $\cA_n$ lives on the space of continuous functions. 
\end{thm}

The local comparison technique can be used to study local fluctuations of last-passage times for lengths of size $n^{\gamma}$, with $\gamma\in(0,2/3)$ (so smaller than the typical scale $n^{2/3}$). Let $\cB\equiv(\cB(u))_{u\geq 0}$ denote the standard two-sided Brownian motion process.
\begin{thm}\label{thm:LocalFluct}
Fix $\gamma\in(0,2/3)$ and $s>0$ and define $\Delta_n$ by   
$$u\in\R\,\mapsto\, \Delta_n(u):=\frac{L[sn+un^{\gamma}]_{n}-L[sn]_n-\mu un^{\gamma}}{\sigma n^{\gamma/2}}\,,$$
where $\mu:=s^{-1/2}$ and $\sigma:=s^{-1/4}$. Then 
$$\lim_{n\to\infty}\Delta_n(\cdot)\stackrel{dist.}{=}\cB(\cdot)\,,$$
in the sense of weak convergence of probability measures in the space of cadlag functions.
\end{thm}

As we mentioned in the previous section, the Airy process locally behaves like Brownian motion \cite{CH,H}. By applying the local comparison technique again, we will present an alternative proof of the functional limit theorem for this local behavior. 
\begin{thm}\label{thm:LocalAiry}
Define $\cA^\epsilon$ by  
$$u\in\R\,\mapsto\,\cA^{\epsilon}(u):=\epsilon^{-1/2}\left(\cA(\epsilon u)-\cA(0)\right)\,.$$
Then
$$\lim_{\epsilon\to 0}\cA^{\epsilon}(\cdot)\stackrel{dist.}{=}\sqrt{2}\cB(\cdot)\,,$$ 
in the sense of weak convergence of probability measures in the space of continuous functions.
\end{thm}
          
\subsection{The lattice model with exponential weights}          
Consider a collection $\{\omega_{[x]_t}\,:\,[x]_t\in\ZZ^2\}$ of  i.i.d. non negative random variables with an exponential distribution of parameter one.  Let $\Pi([x]_t,[y]_u)$ denote the collection of all lattice paths $\varpi=([z]_{v_j})_{j=1,\dots,k}$ such that: 
\begin{itemize}
\item $[z]_{v_1}\in\{[x]_t+[1]_0,[x]_t+[0]_1\}$ and $[z]_{v_k}=[y]_u$;
\item $[z]_{v_{j+1}}-[z]_{v_j}\in\{[1]_0,[0]_1\}$ for $j=0,1\dots,,k-1$. 
\end{itemize}
The (lattice) last-passage percolation time between $[x]_t <[y]_u$ is defined by 
$$L^l([x]_t,[y]_u):=\max_{\varpi\in\Pi([x]_t,[y]_u)}\big\{\sum_{[z]_v\in\varpi}\omega_{[z]_v}\big\}\,.$$
Denote $L^l[x]_t:=L^l([0]_0,[x]_t)$ and define $\cA^l_n$ by 
$$u\in\R\,\mapsto\,\cA^l_n(u):=\frac{L^l[n+2^{5/3}un^{2/3}]_n-(4n+2^{8/3}un^{2/3})+2^{4/3}u^2n^{1/3}}{2^{4/3}n^{1/3}}\,.$$
Corwin, Ferrari and P\'ech\'e \cite{CFP} proved that 
\begin{equation}\label{eq:Airy}
\lim_{n\to\infty}\cA^l_n(\cdot)\stackrel{dist.}{=}\cA(\cdot)\,,
\end{equation}
in the sense of finite dimensional distributions. The local comparison method can be used in this context as well. (The lattice version of Lemma \ref{lem:LocalComparison} is straightforward. For exponential weights, the analog to Lemma \ref{lem:ExitControl} was proved in \cite{BCS}.) 

\begin{thm}\label{thm:LTight}
The collection $\{\cA^l_n\}$ is tight in the space of cadlag functions on $[a,b]$. Furthermore, any weak limit of $\cA_n$ lives on the space of continuous functions. 
\end{thm}

\begin{thm}\label{thm:LLocalFluct}
Fix $\gamma\in(0,2/3)$ and $s>0$ and define $\Delta_n$ by   
$$u\in\R\,\mapsto\, \Delta^l_n(u):=\frac{L^l[sn+un^{\gamma}]_{n}-L^l[sn]_n-\mu un^{\gamma}}{\sigma n^{\gamma/2}}\,,$$
where $\mu=\sigma:=s^{-1/2}(1+s^{1/2})$. Then 
$$\lim_{n\to\infty}\Delta^l_n(\cdot)\stackrel{dist.}{=}\cB(\cdot)\,,$$
in the sense of weak convergence of probability measures in the space of cadlag functions.
\end{thm}          
 
To avoid repetitions, we will not present a proof of the lattice results. We hope that the reader can convince his (or her) self that the method that we will describe in detail for the Hammersley last-passage percolation model can be easily adapted to the lattice models with exponentials weights. 
            
\begin{rem}
We also expect that the local comparison method can be used in the log-gamma polymer model, introduced by Sepp\"al\"ainen \cite{T}. The polymer versions of Lemma  \ref{lem:LocalComparison} and Lemma \ref{lem:ExitControl} were proved in \cite{T}.  
\end{rem}
               
\section{Local comparison and exit points}     

The Hammersley last-passage percolation model has a representation as an interacting particle system, called the Hammersley process  \cite{AD,CG1}. We will use notations used in \cite{CP}. In the positive time half plane we have the same planar Poisson point process as before. On the $x$-axis we take a Poisson process of intensity $\lambda>0$. The two Poisson process are assumed to be independent of each other. For $x\in \R$ and $t>0$ we define
$$L_{\la}[x]_t\equiv L_{\nu_\lambda}[x]_t:=\sup_{z\in (-\infty,x]} \left\{ \nu_\la(z) + L([z]_0,[x]_t)\right\}\,,$$
where, for $z\leq x$,
$$\nu_\la(z)=\left\{\begin{array}{ll}\mbox{ the number of Poisson points in }(0,z]\times\{0\}& \mbox{ for } z> 0\\
\mbox{ minus the number of Poisson points in }(z,0]\times \{0\} & \mbox{ for } z\leq0\,.\end{array}\right.$$ 
The process $(M^t_\lambda)_{t\geq 0}$, given by
$$M^t_{\nu_\la}(x,y]\equiv M^t_\la(x,y]:=L_\la[y]_t-L_\la[x]_t\,\,\,\mbox{ for }x<y\,,$$
is a Markov process on the space of locally finite counting measures on $\R$. The Poisson process is the invariant measure of this particle system in the sense that
\begin{equation}\label{eq:equilibrium} 
M^t_\la\stackrel{dist.}{=}\mbox{ Poisson process of intensity $\lambda$ for all $t\geq 0$}\,.
\end{equation}
Notice that the last-passage time $L=L_{\nu_0}$ can be recovered in the positive quadrant by choosing a measure $\nu_0$ on the axis that has no points to the right of $0$, and an infinite amount of points in every interval $(-\eps,0)$, $\forall \eps>0$ (this could be called a ``wall'' of points). Thus 
$$L_\la\equiv L_{\nu_\lambda}({\cal P})\,\mbox{ and }\, L\equiv L_{\nu_0}({\cal P})\,$$
are coupled by the same two-dimensional Poisson point process $\cal P$, which corresponds to the basic coupling between $M^t_{\nu_\lambda}$ and $M^t_{\nu_0}$. 

Define the exit points 
$$Z_\la[x]_t :=\sup\left\{z\in(-\infty,x]\,:\,L_\la[x]_t=\nu_\lambda(z)+L([z]_0,[x]_t)\right\}\,,$$
and
$$Z'_\la[x]_t :=\inf\left\{z\in(-\infty,x]\,:\,L_\la[x]_t=\nu_\lambda(z)+L([z]_0,[x]_t)\right\}\,.$$
By using translation invariance and invariance under the map $(x,t)\mapsto(\lambda x,t/\lambda)$, we have that 
\begin{equation}\label{eq:sym}  
Z_\lambda[x+h]_t\stackrel{dist.}{=}Z_\lambda[x]_t+h\,\,\,\mbox{ and }\,\,\, Z_\la[x]_t\stackrel{dist.}{=}\lambda Z_1[\lambda x]_{t/\lambda}.
\end{equation}
We need to use one more symmetry. In \cite{CG1}, the Hammersley process was set up as a process in the first quadrant with sources on the $x$-axis and sinks on the $t$-axis. In our notation this means that the process $t\mapsto L_\lambda[0]_t$ is a Poisson process with intensity $1/\lambda$ which is independent of the Poisson process $\nu_\lambda$ restricted to the positive $x$-axis and independent of the Poisson process in the positive quadrant. We can now use reflection in the diagonal to see that the following equality holds:
\begin{equation}\label{eq:diagsym}
\P\left(Z'_\lambda[x]_t < 0\right) = \P\left(Z_{1/\lambda}[t]_x > 0\right).
\end{equation}
We use that $Z'_\lambda[x]_t<0$ is equivalent to the fact that the maximizing path to the left-most exit point crosses the positive $t$-axis, and not the positive $x$-axis. 

The local comparison technique consists of bounding from above and from below the local differences of $L$ by the local differences of $L_\lambda$. These bounds depend on the position of the exit points. It is precisely summarized by the following lemma.  
\begin{lem}\label{lem:LocalComparison}
Let $0\leq x\leq y$ and $t\geq 0$. If $Z'_\la[x]_t\geq0$ then
$$L[y]_t-L[x]_t\leq L_\la[y]_t-L_\la[x]_t\,,$$
and if $Z_\la[y]_t\leq 0$ then 
$$L[y]_t-L[x]_t\geq L_\la[y]_t-L_\la[x]_t\,.$$
\end{lem}

\noindent{\bf Proof\,\,} When we consider a path $\varpi$ from $[x]_s$ to $[y]_t$ consisting of increasing points, we will view $\varpi$ as the lowest increasing continuous path connecting all the points, starting at $[x]_s$ and ending at $[y]_t$. In this way we can talk about crossings with other paths or with lines. The geodesic between $[x]_s$ and $[y]_t$ is given by the lowest path (in the sense we just described) that attains the maximum in the definition of $L([x]_s,[y]_t)$. We will denote this geodesic by $\varpi([x]_s,[y]_t)$. Notice that
$$L([x]_s,[y]_t)=L([x]_s,[z]_r)+L([z]_r,[y]_t)\,,$$
for any $[z]_r\in\varpi([x]_s,[y]_t)$. 

Assume that $Z'_\la[x]_t\geq 0$ and let $\cc$ be a crossing between the two geodesics $\varpi([0]_0,[y]_t)$ and $\varpi([z']_0,[x]_t)$, where $z':=Z'_\la[x]_t$. Such a crossing always exists because $x\leq y$ and $z'=Z'_\la[x]_t\geq 0$.  We remark that, by superaddivity,
$$L_\la[y]_t  \geq  \nu_\la(z') + L([z']_0,[y]_t) \geq  \nu_\la(z') + L([z']_0,\cc) + L(\cc,[y]_t)\,.$$
We use this, and that (since $\cc\in\varpi([z']_0,[x]_t)$)
$$ \nu_\la(z') + L([z']_0,\cc)-L_\la [x]_t= -L(\cc,[x]_t)\,,$$
in the following inequality:
\begin{eqnarray*}
L_\la[y]_t - L_\la[x]_t & \geq & \nu_\la\big(z'\big)+L([z']_0,\cc) + L(\cc,[y]_t) - L_\la[x]_t\\
& = & L(\cc,[y]_t) - L(\cc,[x]_t)\,.
\end{eqnarray*}
By superaddivity,
$$ - L(\cc\,,\,[x]_t)\geq L([0]_0,\cc)-L[x]_t\,,$$
and hence (since $\cc\in\varpi([0]_0,[y]_t)$)
\begin{eqnarray*}
L_\la[y]_t - L_\la[x]_t & \geq & L(\cc,[y]_t) - L(\cc,[x]_t)\\
& \geq & L(\cc,[y]_t) + L([0]_0,\cc)-L([0]_0,[x]_t)\\
& = & L[y]_t-L[x]_t\,.
\end{eqnarray*}
The proof of the second inequality is very similar. Indeed, denote $z:=Z_\la[y]_t$ and let $\cc$ be a crossing between $\varpi([0]_0,[x]_t)$ and $\varpi([z]_0,[y]_t)$. By superaddivity,
$$L_\la[x]_t \geq \nu_\la(z) + L([z]_0,[x]_t) \geq  \nu_\la(z) + L([z]_0,\cc) + L(\cc,[x]_t)\,.$$
Since $\cc\in\varpi([z]_0,[y]_t)$ we have that 
$$L_\la[y]_t-\nu_\la(z) - L([z]_0,\cc)=L(\cc,[y]_t)\,,$$
which implies that 
\begin{eqnarray*}
L_\la[y]_t - L_\la[x]_t & \leq &L_\la[y]_t- \nu_\la(z)-L([z]_0,\cc) - L(\cc,[x]_t)\\
& = & L(\cc,[y]_t) - L(\cc,[x]_t)\\
& \leq & L[y]_t-L([0]_0,\cc)-  L(\cc,[x]_t)\\
& = & L[y]_t-L[x]_t\,,
\end{eqnarray*}
where we have used that $\cc\in\varpi([0]_0,[x]_t)$ in the last step.

\hfill$\Box$\\ 

\begin{rem} In fact the first statement of the lemma is also true when $Z_\lambda[x]_t\geq 0$ and the second statement is true when $Z'_\lambda[y]_t\leq 0$, both without any change to the given proof. This is a stronger statement, since $Z'_\lambda[x]_t\leq Z_\lambda[x]_t$, but we will only need the lemma as it is formulated.
\end{rem}

In order to apply Lemma \ref{lem:LocalComparison} and extract good bounds for the local differences one needs to control the position of exit points. This is given  by the next lemma.
\begin{lem}\label{lem:ExitControl}
There exist constant $C>0$ such that,
$$\P\left(Z_1[n]_n > r n^{2/3}\right)\leq \frac{C}{r^3}\,,$$
for all $r\geq 1$ and all $n\geq 1$.
\end{lem}

\noindent{\bf Proof\,\,} See Corollary 4.4 in \cite{CG2}.

\hfill$\Box$\\ 

\section{Proof of Theorem \ref{thm:Tight}}   
For simple notation, and without loss of generality, we will restrict our proof to $[a,b]=[0,1]$.
\begin{lem}\label{lem:tight}
Fix $\beta\in(1/3,1)$ and for each $\delta\in(0,1)$ and $n\geq 1$ set
$$\lambda_\pm=\lambda_\pm(n,\delta):=1\pm\frac{\delta^{-\beta}}{n^{1/3}}\,.$$
Define the event  
$$E_n(\delta):=\left\{Z'_{\lambda_{+}}[n]_n\geq 0\,\,\mbox{ and }\,\,Z_{\lambda_{-}}[n+2n^{2/3}]_n\leq 0\right\}\,.$$
Then
there exists a constant $C>0$ such that, for sufficiently small $\delta>0$,
$$\limsup_{n\to\infty}\P\left(E_n(\delta)^c\right)\leq C\delta^{3\beta}\,.$$
\end{lem}

\noindent{\bf Proof\,\,} Denote $r:=\delta^{-\beta}$ and let 
$$n_+:=\lambda_+ n< 2n\,\,\mbox{ and }\,\,h_{+}:=\left(\lambda_+-\frac{1}{\lambda_+}\right) n> rn^{2/3}> r n_+^{2/3}/2 \,$$
(for all sufficiently large $n$). By \eqref{eq:sym} and \eqref{eq:diagsym},
\begin{eqnarray*}
\P\left(Z'_{\lambda_+}[n]_n<0\right)&=&\P\left(Z_{1/\lambda_+}[n]_n>0\right)\\
&=&\P\left(Z_1[n/\lambda_+]_{\lambda_+ n}>0\right)\\
&=&\P\left(Z_1[\lambda_+n-h_+]_{\lambda_+ n}>0\right)\\
&=&\P\left(Z_1[\lambda_+n]_{\lambda_+ n}>h_+\right)\\
&\leq &\P\left(Z_1[n_+]_{n_+}>r n_+^{2/3}/2\right)\,.
\end{eqnarray*}
Analogously, for
$$n_-:= \frac{n}{\lambda_-}<2n\,\,\mbox{ and }\,\,h_{-}:=\left(\frac{1}{\lambda_-}-\lambda_-\right) n> rn^{2/3}>r n_-^{2/3}/2 \,,$$
we have that 
\begin{eqnarray*}
\P\left(Z_{\lambda_-}[n+2 n^{2/3}]_n>0\right)&=&\P\left(Z_{\lambda_-}[n]_n>-2 n^{2/3}\right)\\
&=&\P\left(\lambda_-Z_{1}[\lambda_-n]_{n/\lambda_-}>-2 n^{2/3}\right)\\
&\leq&\P\left(Z_{1}[n_- - h_-]_{n/\lambda_-}>-2 n^{2/3}\right)\\
&=&\P\left(Z_{1}[n_-]_{n_-}>h_-  -2 n^{2/3}\right)\\
&\leq& \P\left(Z_1[n_-]_{n_-}>(r-4) n_-^{2/3}/2\right)\,.
\end{eqnarray*}
Now one can use Lemma \ref{lem:ExitControl} to finish the proof.

\hfill$\Box$\\

\begin{lem}\label{lem:compa}
Let $\delta\in(0,1)$ and $u\in[0,1-\delta)$. Then, on the event $E_n(\delta)$, for all $v\in[u,u+\delta]$ we have that
$$ \cB_{n,-}(v)-\cB_{n,-}(u)-2\delta^{1-\beta}\leq\cA_n(v)-\cA_n(u)\,\leq\, \cB_{n,+}(v)-\cB_{n,+}(u)+4\delta^{1-\beta}\,,$$
where 
$$\cB_{n,\pm}(u):\frac{L_{\lambda\pm}[n+2un^{2/3}]_n-L_{\lambda\pm}[n]_n-\lambda_\pm 2un^{2/3}}{n^{1/3}}\,.$$
\end{lem}

\noindent{\bf Proof\,\,} For fixed $t$, $Z_\lambda'[x]_t$ and $Z_\lambda[x]_t$ are non-decreasing functions of $x$. Thus,  on the event $E_n(\delta)$, 
$$Z'_{\lambda_+}[n+2un^{2/3}]_n\geq 0\,\mbox{ and }\,Z_{\lambda_-}[n+2(u+\delta)n^{2/3}]_n\leq 0\,.$$
By Lemma \ref{lem:LocalComparison}, this implies that, for all $v\in[u,u+\delta]$,  
$$L[n+vn^{2/3}]_{n}-L[n+un^{2/3}]_n\leq L_{\lambda_+}[n+vn^{2/3}]_{n}-L_{\lambda_+}[n+un^{2/3}]_n\,,$$
and
$$L[n+vn^{2/3}]_{n}-L[n+un^{2/3}]_n\geq L_{\lambda_-}[n+vn^{2/3}]_{n}-L_{\lambda_-}[n+un^{2/3}]_n\,.$$
Since 
$$(\lambda_+-1)(2v-2u)n^{1/3}+v^2-u^2\leq 2\delta^{1-\beta}+2\delta\leq 4\delta^{1-\beta}\,,$$
and
$$(\lambda_- -1)(2v-2u)n^{1/3}+v^2-u^2\geq -2\delta^{1-\beta}\,,$$ 
we have that, on the event  $E_n(\delta)$,
$$\cA_n(v)-\cA_n(u)\,\leq\, \cB_{n,+}(v)-\cB_{n,+}(u)+4\delta^{1-\beta}\,.$$
and
$$\cA_n(v)-\cA_n(u)\,\geq\, \cB_{n,-}(v)-\cB_{n,-}(u)-2\delta^{1-\beta}\,,$$
for all $v\in[u,u+\delta]$.

\hfill$\Box$\\

\noindent{\bf Proof of Theorem \ref{thm:Tight}\,\,}
For fixed $u\in[0,1)$ take $\delta>0$ such that $u+\delta\leq1$. By Lemma \ref{lem:compa},
$$\sup_{v\in[u,u+\delta]}|\cA_n(v)-\cA_n(u)|\leq \max\left\{\sup_{v\in[u,u+\delta]}|\cB_{n,\pm}(v)-\cB_{n,\pm}(u)| \right\}+4\delta^{1-\beta}\,,$$
on the event $E_n(\delta)$. Hence, for any $\eta>0$,
\begin{eqnarray*}
\P\left( \sup_{v\in[u,u+\delta]}|\cA_n(v)-\cA_n(u)|>\eta\right)&\leq&\P\left(E_n(\delta)^c\right)\\
&+&\P\left( \sup_{v\in[u,u+\delta]}|\cB_{n,+}(v)-\cB_{n,+}(u)|>\eta-4\delta^{1-\beta}\right)\\
&+&\P\left( \sup_{v\in[u,u+\delta]}|\cB_{n,-}(v)-\cB_{n,-}(u)|>\eta-4\delta^{1-\beta}\right)\,.
\end{eqnarray*}
By \eqref{eq:equilibrium},
$$P_n(x):=L_\lambda\big([n+x]_n\big)-L_\lambda\big([n]_n\big)\,,\,\mbox{ for }\,\,\,x\geq 0\,,$$
is a Poisson process of intensity $\lambda$. Since $\lambda^{\pm}\to1$ as $n\to\infty$, $\cB_{n,-}(u/2)$ and $\cB_{n,+}(u/2)$ converge in distribution to a standard Brownian motion $\cB$. Thus,  by Lemma \ref{lem:tight}, for $\delta<(\eta/8)^{1/(1-\beta)}$,
\begin{eqnarray*}
\limsup_{n\to\infty}\P\left(\sup_{v\in[u,u+\delta]}|\cA_n(v)-\cA_n(u)|>\eta\right)&\leq& C\delta^{3\beta}+2\P\left( \sup_{v\in[u,u+\delta]}|\cB(2v)-\cB(2u)|>\eta-4\delta^{1-\beta}\right)\\
&\leq&C\delta^{3\beta}+2\P\left( \sup_{v\in[0,1]}|\cB(v)|>\frac{\eta}{2\sqrt{2\delta}}\right)\,,
\end{eqnarray*}
which implies that (recall that $\beta\in(1/3,1)$)
\begin{equation}\label{tight}
\limsup_{\delta\to 0^+}\frac{1}{\delta}\left(\limsup_{n\to\infty}\P\left(\sup_{v\in[u,u+\delta]}|\cA_n(v)-\cA_n(u)|>\eta\right)\right)=0\,.
\end{equation}
Since \cite{CG2}
$$\limsup_{n\to\infty}\frac{\E |L[n]_n-2n|}{n^{1/3}}<\infty\,,$$
we have that $\{\cA_n(0)\,,\,n\geq 1\}$ is tight. Together with \eqref{tight}, this shows tightness of the collection $\{\cA_n\,,\,n\geq 1\}$ in the space of cadlag functions on $[0,1]$, and also that every weak limit lives in the space of continuous functions \cite{Bi}. 

\hfill$\Box$\\ 

\section{Proof of Theorem \ref{thm:LocalFluct}} 
For simple notation, we will prove the statement for $s=1$ and restrict our selves to $[0,1]$. The reader can then check that rescaling gives the result for general $s>0$, since
\[ L[sn]_n \stackrel{dist.}{=} L[s^{1/2}n]_{s^{1/2}n}.\]
\begin{lem}\label{lem:PoissonTight}
Fix $\gamma'\in(\gamma,2/3)$ and let 
$$\la_{\pm}=\la_{\pm}(n):=1\pm\frac{1}{n^{\gamma'/2}}\,.$$
Define the event
$$E_n:=\left\{Z'_{\la_+}[n]_n\geq0\,\,\mbox{ and }\,\,Z_{\la_-}[n+n^{\gamma}]_n\leq0\right\}\,.$$
There exists a constant $C>0$ such that 
$$\P\left(E_n^c\right)\leq \frac{C}{n^{1-3\gamma'/2}}\,$$
for all sufficiently large $n$.
\end{lem}

\noindent{\bf Proof\,\,} Denote $r:=n^{1/3-\gamma'/2}$ and let 
$$n_+:=\lambda_+ n< 2n\,\,\mbox{ and }\,\,h_{+}:=\left(\lambda_+-\frac{1}{\lambda_+}\right) n> rn^{2/3}> r n_+^{2/3}/2 \,,$$
(for all sufficiently large $n$). By \eqref{eq:sym} and \eqref{eq:diagsym},
\begin{eqnarray*}
\P\left(Z'_{\lambda_+}[n]_n<0\right)&=&\P\left(Z_{1/\lambda_+}[n]_n>0\right)\\
&=&\P\left(Z_1[n/\lambda_+]_{\lambda_+ n}>0\right)\\
&=&\P\left(Z_1[\lambda_+n-h_n]_{\lambda_+ n}>0\right)\\
&=&\P\left(Z_1[\lambda_+n]_{\lambda_+ n}>h_+\right)\\
&\leq&\P\left(Z_1[n_+]_{n_+}>r n_+^{2/3}/2\right)\,.
\end{eqnarray*}
Analogously, for
$$n_-:= \frac{n}{\lambda_-}<2n\,\,\mbox{ and }\,\,h_{-}:=\left(\frac{1}{\lambda_-}-\lambda_-\right) n> rn^{2/3}> r n_-^{2/3}/2 \,,$$
we have that 
\begin{eqnarray*}
\P\left(Z_{\lambda_-}[n+n^{\gamma}]_n>0\right)&=&\P\left(Z_{\lambda_-}[n]_n>-n^{\gamma}\right)\\
&=&\P\left(\lambda_-Z_{1}[\lambda_-n]_{n/\lambda_-}>-n^{\gamma}\right)\\
&\leq&\P\left(Z_{1}[n_-]_{n_-}>h_- - n^{\gamma}\right)\\
&\leq& \P\left(Z_1[n_-]_{n_-}>(r-n^{\gamma-2/3}) n_-^{2/3}/2\right)\,,
\end{eqnarray*}
Now one can use Lemma \ref{lem:ExitControl} to finish the proof.

\hfill$\Box$\\ 

\begin{lem}\label{lem:PoissonComparison}
On the event $E_n$, for all $u<v$ in $[0,1]$,  
$$\Gamma_n^{-}(v)-\Gamma_n^{-}(u)-\frac{1}{n^{(\gamma'-\gamma)/2}}\leq \Delta_n(v)-\Delta_n(u)\leq \Gamma_n^+(v)-\Gamma_n^+(u)+\frac{1}{n^{(\gamma'-\gamma)/2}}\,,$$
where 
$$\Gamma_n^{\pm}(u):=\frac{L_{\la_\pm}[n+un^{\gamma}]_{n}-L_{\la_\pm}[n]_n-\lambda_\pm un^{\gamma}}{n^{\gamma/2}}\,.$$

\end{lem}

\noindent{\bf Proof\,\,} By Lemma \ref{lem:LocalComparison}, if $Z'_{\lambda_+}[n]_n\geq 0$ then  
$$L[n+vn^{\gamma}]_{n}-L[n+un^{\gamma}]_n\leq L_{\lambda_+}[n+vn^{\gamma}]_{n}-L_{\lambda_+}[n+un^{\gamma}]_n\,,$$
and if $Z_{\lambda_-}[n+n^{\gamma}]_n\leq 0$ then
$$L[n+vn^{\gamma}]_{n}-L[n+un^{\gamma}]_n\geq L_{\lambda_-}[n+vn^{\gamma}]_{n}-L_{\lambda_-}[n+un^{\gamma}]_n\,.$$
Using that $\lambda_{\pm}:=1\pm n^{-\gamma'/2}$, one can finish the proof of the lemma. 

\hfill$\Box$\\ 

\noindent{\bf Proof of Theorem \ref{thm:LocalFluct}\,\,} By Lemma \ref{lem:PoissonComparison}, on the event $E_n^c$,
$$| \Delta_n(v)-\Delta_n(u)|\leq\max \left\{|\Gamma_n^{\pm}(v)-\Gamma_n^{\pm}(u)|\right\}+\frac{1}{n^{(\gamma'-\gamma)/2}}\,.$$
Thus, by Lemma \ref{lem:PoissonTight},
\begin{eqnarray*}
\P\left(\sup_{v\in[u,u+\delta]} | \Delta_n(v)-\Delta_n(u)|>\eta\right)&\leq&\P\left(\sup_{v\in[u,u+\delta]} |\Gamma^+_n(v)-\Gamma^+_n(u)|+\frac{1}{n^{(\gamma'-\gamma)/2}}>\eta\right)\\
&+&\P\left(\sup_{v\in[u,u+\delta]} |\Gamma^-_n(v)-\Gamma^-_n(u)|+\frac{1}{n^{(\gamma'-\gamma)/2}}>\eta\right)\\
&+&\P\left(E_n^c\right)\,.
\end{eqnarray*}
As before, $\lambda^{\pm}\to1$ as $n\to\infty$, which implies that 
$$\limsup_{n\to\infty}\P\left(\sup_{v\in[u,u+\delta]} | \Delta_n(v)-\Delta_n(u)|>\eta\right)\leq 2\P\left(\sup_{v\in[0,\delta]} |B(v)|>\eta\right)=2\P\left(\sup_{v\in[0,1]} |B(v)|>\frac{\eta}{\sqrt{\delta}}\right)\,,$$
and hence
\begin{equation}\label{eq:tight}
\limsup_{\delta\to0^+}\frac{1}{\delta}\left(\limsup_{n\to\infty}\P\left(\sup_{v\in[u,u+\delta]} | \Delta_n(v)-\Delta_n(u)|>\eta\right)\right)=0\,.
\end{equation}
Since $\Delta_n(0)=0$, \eqref{eq:tight} implies tightness of the collection $\{\Delta_n\,,\,n\geq 1\}$ in the space of cadlag functions on $[0,1]$, and also that every weak limit lives in the space of continuous functions \cite{Bi}. 

The finite dimensional distributions of the limiting process can be obtained in the same way. Indeed, by Lemma \ref{lem:PoissonTight} and Lemma \ref{lem:PoissonComparison}, for $u_1,\dots,u_k\in[0,1]$ and $a_1,\dots,a_k,\in\R$, 
$$\P\left(\cap_{i=1}^k\left\{\Delta_n(u_i)\leq a_i\right\}\right)\geq \P\left(\cap_{i=1}^k\left\{\Gamma^+_n(u_i)\leq a_i -\frac{1}{n^{(\gamma'-\gamma)/2}}\right\}\right)-\P\left(E_n^c\right)\,,$$
and 
$$\P\left(\cap_{i=1}^k\left\{\Delta_n(u_i)\leq a_i\right\}\right)\leq\P\left(\cap_{i=1}^k\left\{\Gamma^-_n(u_i)\leq a_i+\frac{1}{n^{(\gamma'-\gamma)/2}}\right\}\right)+\P\left(E_n^c\right)\,,$$
which shows that the finite dimensional distributions of $\Delta_n$ converge to the finite dimensional distributions of the standard Brownian motion process. 

\hfill$\Box$\\ 

\section{Proof of Theorem \ref{thm:LocalAiry}}      

\begin{lem}\label{lem:local}
Fix $\beta\in(0,1/2)$ and for $\epsilon\in(0,1)$ let
$$\lambda_\pm=\lambda_\pm(n,\epsilon):=1\pm\frac{\epsilon^{-\beta}}{n^{1/3}}\,.$$
Define the event 
 $$E_n(\epsilon):=\left\{Z'_{\lambda_{+}}([n]_n)\geq 0\,\,\mbox{ and }\,\,Z_{\lambda_{-}}([n+n^{2/3}]_n)\leq 0\right\}\,.$$
There exists a constant $C>0$ such that, for all sufficiently small $\epsilon>0$,
$$\limsup_{n\to\infty}\P\left(E_n(\epsilon)^c\right)\leq C\epsilon^{3\beta}\,.$$
\end{lem}

\noindent{\bf Proof\,\,} The same proof as in Lemma \ref{lem:tight} applies.

\hfill$\Box$\\ 

\begin{lem}\label{lem:localcompa} 
On the event $E_n(\epsilon)$, for all $u\in[0,1-\delta)$ and  $v\in[u,u+\delta]$, we have that
$$ \cB_{n,-}(\epsilon v)-\cB_{n,-}(\epsilon u)-2\delta\epsilon^{1-\beta}\leq\cA_n(\epsilon v)-\cA_n(\epsilon u)\,\leq\, \cB_{n,+}(\epsilon v)-\cB_{n,+}(\epsilon u)+4\delta\epsilon^{1-\beta}\,.$$
\end{lem}

\noindent{\bf Proof\,\,} The same proof as in Lemma \ref{lem:compa} applies. Note that in this case we have
\[ (\lambda_+ - 1)(2\eps v- 2\eps u) \leq 2\eps^{1-\beta}\delta.\]

\hfill$\Box$\\

\noindent{\bf Proof of Theorem \ref{thm:LocalAiry}\,\,}
For $u\in[0,1]$, let 
$$\cA_n^\epsilon(u):=\epsilon^{-1/2}\left(\cA_n(\epsilon u)-\cA_n(0)\right)\,\mbox{ and }\,\cB^\epsilon_{n,\pm}(u):=\epsilon^{-1/2}\cB_{n,\pm}(\epsilon u)\,.$$
By Lemma \ref{lem:localcompa}, on the event $E_n(\epsilon)$, for all $v\in[u,u+\delta]$,
$$ \cB^\epsilon_{n,-}(v)-\cB^\epsilon_{n,-}(u)-2\delta\epsilon^{1/2-\beta}\leq\cA^\epsilon_n(v)-\cA^\epsilon_n(u)\,\leq\, \cB^\epsilon_{n,+}(v)-\cB^\epsilon_{n,+}(u)+4\delta\epsilon^{1/2-\beta}\,,$$
which shows that 
$$\sup_{v\in[u,u+\delta]}|\cA^\epsilon_n(v)-\cA^\epsilon_n(u)|\leq \max\left\{\sup_{v\in[u,u+\delta]}|\cB_{n,\pm}^\epsilon(v)-\cB_{n,\pm}^\epsilon(u)| \right\}+4\delta\epsilon^{1/2-\beta}\,.$$
Therefore,
\begin{eqnarray*}
\P\left( \sup_{v\in[u,u+\delta]}|\cA^\epsilon_n(v)-\cA^\epsilon_n(u)|>\eta\right)&\leq&\P\left(E_n(\epsilon)^c\right)\\
&+&\P\left( \sup_{v\in[u,u+\delta]}|\cB_{n,+}^\epsilon(v)-\cB_{n,+}^\epsilon(u)|>\eta-4\delta\epsilon^{1/2-\beta}\right)\\
&+&\P\left( \sup_{v\in[u,u+\delta]}|\cB_{n,-}^\epsilon(v)-\cB_{n,-}^\epsilon(u)|>\eta-4\delta\epsilon^{1/2-\beta}\right)\,.
\end{eqnarray*}
Since $\cA_n^\epsilon$ is converging to $\cal A^\epsilon$, and $\cB_{n,\pm}^\epsilon$ is converging to a Brownian motion $\cB$, the preceding inequality implies that  
$$\P\left(\sup_{v\in[u,u+\delta]}|\cA^\epsilon(v)-\cA^\epsilon(u)|>\eta\right)\leq C\epsilon^{3\beta}+2\P\left( \sup_{v\in[u,u+\delta]}|\cB(2v)-\cB(2u)|>\eta-4\delta\epsilon^{1/2-\beta}\right)\,.$$
Hence
$$\limsup_{\epsilon\to 0^+}\P\left( \sup_{v\in[u,u+\delta]}|\cA^\epsilon(v)-\cA^\epsilon(u)|>\eta\right)\leq 2\P\left( \sup_{v\in[0,1]}|\cB(v)|>\frac{\eta}{\sqrt{2\delta}}\right)\,,$$
which shows that
\begin{equation}\label{eq:BlowTight}
\limsup_{\delta\to 0^+}\frac{1}{\delta}\left(\limsup_{\epsilon\to 0^+}\P\left( \sup_{v\in[u,u+\delta]}|\cA^\epsilon(v)-\cA^\epsilon(u)|>\eta\right)\right)=0\,.
\end{equation}
Since $\cA^{\epsilon}(0)=0$, by \eqref{eq:BlowTight} we have that $\left\{\cA^\epsilon\,,\,\epsilon\in(0,1]\right\}$ is tight \cite{Bi}.

The finite dimensional distributions of the limiting process can be obtained in the same way. Indeed, by Lemma \ref{lem:localcompa}, for $u_1,\dots,u_k\in[0,1]$ and $a_1,\dots,a_k,\in\R$,
$$\P\left(\cap_{i=1}^k\left\{\cA^\epsilon_n(u_i)\leq a_i\right\}\right)\leq\P\left(\cap_{i=1}^k\left\{\cB^\epsilon_{-,n}(u_i)\leq a_i+4\epsilon^{1/2-\beta}\right\}\right)+\P\left(E_n(\epsilon)^c\right)\,,$$
and
$$\P\left(\cap_{i=1}^k\left\{\cA^\epsilon_n(u_i)\leq a_i\right\}\right)\geq \P\left(\cap_{i=1}^k\left\{\cB^\epsilon_{+,n}(u_i)\leq a_i-4\epsilon^{1/2-\beta}\right\}\right)-\P\left(E_n(\epsilon)^c\right)\,.$$
Thus, by Lemma \ref{lem:local}, 
$$ \P\left(\cap_{i=1}^k\left\{\cA^\epsilon(u_i)\leq a_i\right\}\right)\leq \P\left(\cap_{i=1}^k\left\{\cB(2u_i)\leq a_i+4\epsilon^{1/2-\beta}\right\}\right)+C\epsilon^{3\beta}\,,$$
and
$$ \P\left(\cap_{i=1}^k\left\{\cA^\epsilon(u_i)\leq a_i\right\}\right)\geq \P\left(\cap_{i=1}^k\left\{\cB(2u_i)\leq a_i-4\epsilon^{1/2-\beta}\right\}\right)-C\epsilon^{3\beta}\,,$$
which proves that,
\begin{equation}\label{eq:BlowDist} 
\lim_{\epsilon\to 0^+}\P\left(\cap_{i=1}^k\left\{\cA^\epsilon(u_i)\leq a_i\right\}\right)=\P\left(\cap_{i=1}^k\left\{\sqrt{2}\cB(u_i)\leq a_i\right\}\right)\,.
\end{equation}

\hfill$\Box$\\

\end{document}